\documentclass{article}
\def\uu{\bigsqcup}

\title{\bf \Large  Poisson suspensions without roots}
\author{ Valery V. Ryzhikov}
\date{}

\begin{document}

\maketitle


{\normalsize  The discrete rational component in spectrum of an ergodic automorphism S prevents some roots from existing. If S is tensorly multiplied by an ergodic automorphism of the space with a sigma-finite measure, discrete spectrum disappears in this product, but like the smile of Cheshire Cat, the memory of it can remain in the form of the absence of roots. In additional conditions, this effect is inherited by the Poisson suspension over the above product. Starting from this idea, but without using the tensor product, we describe a simple rank-one  construction  for which the Poisson suspension  is rigid and has no roots. }

 
\large
\section{Introduction}
Let $T$ be an automorphism of the space $(X,\mu)$ with sigma-finite measure
$\mu$,
We also denote the corresponding orthogonal operator in $L_2(\mu)$ by $T$.
The group of all automorphisms of the space $(X,\mu)$ embeds continuously into the group of automorphisms of the space of configurations with a probability Poisson measure (see \cite{Ne},\cite{Ro}). The Poisson  image of an automorphism $T$  is called a Poisson suspension and is denoted by $P(T)$.
The suspension $P(T)$ over the automorphism $T$ has a spectral counterpart: the Gaussian automorphism $G(T)$ associated with the action of the orthogonal operator $T$ on the space with  Gaussian probability measure. The automorphisms $G(T)$ and $P(T)$ are spectrally isomorphic, but can have significantly different metric (from the word \it measure \rm) properties, including algebraic properties.

It is well known that the suspensions $P(T)$ and $G(T)$ as operators are isomorphic to the operator
$$1\oplus \bigoplus_{n=1}^\infty T^{\odot n},$$
where $1$ is the identity operator in one-dimensional space, $T^{\odot n}$ is the symmetric $n$-power of the operator $T$.
Gaussian automorphisms have roots; their centralizer is large, since it contains the injective image of the very large  centralizer of the operator $T$ in the group of all orthogonal operators on the real space $L_2(\mu)$.

\bf Sufficient condition for  C(P(T))= P(C(T)).  \rm Let $(\perp)$ denote the following property: \it the spectral measure of an automorphism is mutually singular with its convolution powers. \rm This is equivalent to the fact that the operator $T$ does not have a nonzero intertwining operator with tensor powers $T^{\otimes n}$ as  $n>1$. If the centralizer of an automorphism $T$ in the group of all automorphisms of the space $(X,\mu)$ consists only of $T^n$ and $T$ has the property $(\perp)$, then for the Poisson suspension $P(T)$ its centralizer is trivial, $C(P(T))=\{P(T^n)\}$. This fact is a consequence of Proposition 5.2 of E. Roy's work \cite{Ro}:  $(\perp)$ implies $C(P(T))= P(C(T))$.

\bf Mixing suspensions. \rm If $T$ has no roots and $T$ has the property $(\perp)$, then $P(T)$ also has no roots, this contrasts with the fact that every ergodic Gaussian automorphism $G(T)$ is included in a continuum of non-isomorphic flows. Therefore it has a continuum of roots. For more impressive contrasts, see \cite{PR}. Note also that examples of mixing  suspensions with the trivial centralizer are provided by the works \cite{DR},\cite{RT},\cite{R24}. Recall that $P(T)$ has mixing if $T^n\to 0$ for $n\to\infty$. 
In  \cite{R24} we present mixing  rank one automorphisms with a singular spectrum whose convolution powers are Lebesgue, so, these examples have   the property $(\perp)$ and the trivial centralizer.    

\bf Non-mixing  suspensions \rm $P(T)$ with the trivial centralizer are obtained from the constructions of non-rigid automorphisms proposed in \cite{R19}.

\bf Rigid weakly mixing automorphism without roots. \rm
If the automorphism $T$ is rigid, i.e. $T^{n_i}\to I$ for some
sequences $n_i\to\infty$, the centralizer $C(P(T))$ is large. But if the automorphism $T$ possesses   the property $(\perp)$ and has no roots,   then   $P(T)$ also has no roots.
The typical automorphism of the  probability space is weakly mixing  and  rigid. J. King  \cite{K} proved that it has roots of all orders. So, the constructing a rigid weakly mixing automorphism without roots    requires some ingenuity. The possibility of such constructions by the method of skew products was hinted at by the long-standing work of A.M. Stepin \cite{St} on automorphisms without square roots. T. Adams advised the author that a suitable example was proposed by A. del Junco and D. Rudolph  \cite{JR}. A simple modification of this rank one example ensures that the measure of the phase space is infinite, it is not hard to provide the properties of rigidity and $(\perp)$, but the proof of the absence of roots probably needs the use of subtle analysis. 

\bf Poisson suspensions without roots.  \rm  We offer simple rigid rank-one  constructions  
 and a brief  proof   that  their Poisson suspensions   have no roots. The essence of our  approach is clear from the following statement.

\vspace{2mm}
\bf Theorem 1. \it Let an ergodic automorphism $T$ of a space with sigma-finite measure satisfy the condition $(\perp)$, and its powers $T^n$, $n>0$, have
exactly $n$ invariant ergodic sets. Then $T$ and $P(T)$ have no roots. \rm

\vspace{2mm}

Proof.
   The power $T^n$ has $n$ ergodic sets: $X=\uu_0^{p-1} T^i Y$,
  and the restriction of $T^n$ to $Y$ is an ergodic transformation. Let $R^n=T$. The automorphism $R$ commutes with $T$, so the sets $RT^iY$ are also invariant under $T^n$ and  $T^n$ is ergodic on these sets. Then for some $k$ we have
$$RY=T^kY, \ \ TY=R^nY=T^{kn}Y=Y,$$
that is impossible. If $T$ has the property $(\perp)$, then by Proposition 5.2 \cite{Ro} we get that $P(T)$ has no roots. Theorem is proved.

\section{Examples}
 Let's  provide  rigid automorphisms that satisfy the conditions of Theorem 1.

\it On disjointness of spectrum convolution powers. \rm  
A.B. Katok and A.M. Stepin in \cite{KS} considered the strong convergence of  powers of $T$ to the operator $-I$ on a subspace from $L_2$. This implies the disjointness of spectrum of $T$ with its odd  convolution powers and then  provides  counterexamples to the Kolmogorov conjecture   on the group property of the spectrum of automorphisms.
  V.I. Oseledets in  \cite{O} proposed the use of  non-unitary weak limits $aI$, $0<a<1$. The presence of such limits powers of $T$  forces the mutual singularity of all convolution powers of spectrum  for  $T$, which ensures the fulfillment of the property $(\perp)$. We  use this simple useful observation.
 
 \bf Definition of the  rank one constructions. \rm
We fix a natural number $h_1\geq 1$ (the height of the tower at stage $j=1$), 
a sequence $r_j\to\infty$ (the number of columns into which the tower at stage $j$ is virtually cut) and a sequence of integer vectors (spacer parameters)
$$ \bar s_j=(s_j(1), s_j(2),\dots, s_j(r_j-1),s_j(r_j)).$$
At step $j=1$, the half-interval $E_1$ is specified. Let $j$ be determined at step
  system of disjoint half-intervals
$C_j, TC_j,\dots, T^{h_j-1}C_j,$
and on $C_j, \dots, T^{h_j-2}C_j$
the transformation $T$ is a parallel transfer. 
Such a set of half-intervals is called a tower of stage $j$; 
their union is denoted by $X_j$ and is also called a tower.

Let us represent $C_j$ as a disjunctive union of half-intervals
$C_j^i$, $i=1,2,\dots, {r_j},$ of the same length.
For each $i=1,2,\dots, r_j$, consider the so-called column
$$C_j^i, TC_j^i ,T^2 C_j^i,\dots, T^{h_j-1}C_j^i.$$
To each column with number $i$ we add $s_j(i)$ of 
non-intersecting half-intervals (floors) of length equal to the length of the interval $C_j^i$.
The resulting sets of intervals for fixed $i$,$j$ are called superstructured columns $X_{i,j}$. Note that for a fixed $j$, by construction, the columns $X_{i,j}$ do not intersect. Using parallel transfer of intervals, we now define 
the transformation $T$ so that the columns $X_{i,j}$ have the form
$$C_j^i, TC_j^i, \dots, T^{h_j}C_j^i, T^{h_j+1}C_j^i, \dots, T^{h_j+s_j(i)-1}C_j^ i,$$
   and the upper floors of the columns $X_{i,j}$ ($i<r_j$) 
were mapped by parallel translation to the lower floors of the $T$ transformation
floors of columns $X_{i+1,j}$:
$$T^{h_j+s_j(i)}C_j^i = C_j^{i+1}, \ 0<i<r_j.$$
Putting $C_{j+1}= C^1_j$, we notice that all the indicated floors of the superstructured columns in the new notation have the form
$$C_{j+1}, TC_{j+1}, T^2 C_{j+1},\dots, T^{h_{j+1}-1}C_{j+1},$$
  forming a tower of stage $j+1$ height
  $$ h_{j+1} =h_jr_j +\sum_{i=1}^{r_j}s_j(i).$$

Partial definition of transformation $T$
at stage $j$ is preserved at all subsequent stages.
As a result, we obtain the space $X=\cup_j X_j$
and an invertible transformation $T:X\to X$ preserving
the standard Lebesgue measure $\mu$ on $X$.

\bf Set the parameters. \rm   Let $r'_j\to\infty$, let's put $r_j=jr'_j$.

\bf  Rigidity.  \rm For odd stages $j$ let $s_j(i)=0$.  This provides
$$T^{h_{2j+1}}\to I,\ j\to\infty.$$ \rm

\bf The  property $(\perp)$. \rm For even $j$ we set $s_j(i)=0$ for $0<i\leq r_j/2$ and $s_j(i)=h_j$ for $r_j/2<i\leq r_j$.
 This ensures  the disjointness of spectrum convolution powers due to weak convergence
$$T^{h_{2j}}\to_w I/2,\  j\to\infty. $$

\bf The power $T^n$ has exactly $n$ invariant ergodic sets. \rm Due to the choice of $r_j=jr'_j$ we obtain that for all $j>n$
of height $h_j$ all parameters $s_j(i)$ are multiples of $n$. In this case $T^n$ has exactly $n$ invariant ergodic sets (constructions with a similar property were studied in \cite{R23}).
It is obvious that the  space $X$ can be  represented as a union of non-intersecting towers of heights no less than $h_{n+1}$   and being multiples of $n$.  For the set Y, which appeared in the proof of Theorem 1, we take the union of all floors (of these towers) with numbers equal to $0$ modulo $n$.

So, the rigid  automorphism $T$ satisfies all the conditions of Theorem 1. 
Changing the sequence $r_j'\to\infty$, we can easily obtain a continuum of non-isomorphic 
spectrally suspensions $P(T)$. 

\vspace{2mm}
\bf Theorem 2. \it For our automorphisms  $T$  the Poisson suspensions $P(T)$  are rigid and have  no roots. \rm

\vspace{2mm}
\bf Remark and questions. \rm
  Ergodic Poisson suspensions $P$ and Gaussian automorphisms have the property \cite{R}:
\it  there is an ergodic automorphism $S$ for which
$$ \frac 1 N \sum_{n=1}^N P^{-n}SP^n \ \to \ I, \ n\to\infty.$$
  \rm 

Are there  weakly mixing rigid automorphisms
without this property? Is the map from \cite{JR}  suitable?

Do  generic automorphisms of the probability space have  this property?

\vspace{4mm}
\bf Acknowledgments. \rm The author is grateful to T. Adams, S. Elovatsky and I. Podvigin for comments and questions.

\end{document}